\documentclass[11pt,reqno,nocut]{amsart}

\usepackage[utf8]{inputenc}
\usepackage[T1]{fontenc}
\usepackage{amsmath,amssymb,amsthm,mathtools}

\usepackage{mathrsfs,bm,color}
\usepackage{fancyhdr,lastpage}
\usepackage[shortlabels]{enumitem}
\usepackage{calc}
\usepackage{bbm}
\usepackage{graphicx}
\usepackage{upgreek}
\usepackage{lmodern}
\usepackage[normalem]{ulem}
\usepackage{verbatim}
\usepackage{bbm}
\usepackage{stmaryrd}
\usepackage{amsmath}
\usepackage{amssymb}
\usepackage{dsfont}
\usepackage{amsthm}
\usepackage{hyperref}
\usepackage{array}
\hypersetup{
	colorlinks,
	linkcolor={red!50!black},
	citecolor={red!80!black},
	urlcolor={blue!80!black}
}

\usepackage{xcolor}

\usepackage{graphicx}

\usepackage{amsfonts}%
\usepackage{dsfont}

\allowdisplaybreaks

\usepackage{fullpage}

\numberwithin{equation}{section}

\theoremstyle{plain}
\newtheorem{theo}{Theorem}[section]
\newtheorem{lem}{Lemma}[section]

\theoremstyle{definition}

\newtheorem{defn}{Definition}[section]
\newtheorem{rem}{Remark}[section]

\newcommand{\vect}[1]{\ensuremath{\boldsymbol{\mathbf{#1}}}}

\newcommand\given[1][]{\:#1\vert\:}

\newcommand{\E}{\mathbb{E}}

\newcommand{\sS}{\mathscr{S}}

\newcommand{\N}{\mathbb{N}}

\newcommand{\R}{\mathbb{R}}

\DeclareMathOperator{\pP}{\mathbb{P}}

\newcommand{\dd}{\mathop{}\!\mathrm{d}}

\newcommand{\Var}{\operatorname{Var}}

\renewcommand{\epsilon}{\varepsilon}

\begin{document}
	
	\title{Synchronization of Enriched Spin Glasses}
	
	\author{Franziska Engler}
	\address[Franziska Engler]{Syracuse University}
	\email{fengler@syr.edu}
	\author{Justin Ko}
	\address[Justin Ko]{Syracuse University}
	\email{jko104@syr.edu}
	
	\date{\today}

	\begin{abstract}
		We show that enriching a non-convex multi-species Hamiltonian by a uniform Ruelle probability cascade implies synchronization of the joint overlap array without an additional regularizing perturbation. In particular, we show that any limit point of the overlap array is a deterministic non-decreasing function of a uniform Ruelle probability cascade.
	\end{abstract}
	\maketitle

	\section{Introduction}
	
	In this short note, we show that the enrichment of the free energy with a uniform Ruelle probability cascade (RPC) \cite{Ruelle,ruelle2} plays the role of a macroscopic perturbation that forces synchronization in the limit without the need of an additional perturbation argument. The form of the enrichment is special in the sense that it is rich enough to capture the underlying ultrametric \cite{PUltra} structure of the spin configurations in the limit, making it a natural candidate to synchronize with the configurations. One might expect that the synchronization is only possible for optimal enrichments that match the limiting structure, but this is not the case and it holds for any regular (see Definition~\ref{defn:reg}) enrichment path. 
	
	Fragments of this connection appeared in several recent works \cite{ho_converge,boursier_external_field,chen2026formulafreeenergynonconvex} which proved estimates on conditional overlaps by conditioning on the enrichment. However, each paper used a different representation, which was not the usual form of the enrichment that appeared in earlier versions of this model. One objective is to show that the `cascade conditioning' trick also implies synchronization when enriched by the familiar Ruelle probability cascade with continuous uniform distribution. The main objective is to interpret some consequences of the cascade conditioning. 
	
	Analogous to the inverse transform sampling method, we show that any sub-sequential limit of the off-diagonal entries of the overlap arrays $R_{\ell,\ell'}$ from an unperturbed, but enriched, measure are non-decreasing functions of the same continuous uniform Ruelle probability cascade $U_{\ell,\ell'}$, i.e.
	\[
	(R_{\ell,\ell'}, U_{\ell,\ell'})_{\ell\neq \ell'} \stackrel{L^2}{\to} (F^{-1}_R(U_{\ell,\ell'}) , U_{\ell,\ell'})_{\ell\neq \ell'}
	\]
	where $F_R^{-1}$ is the quantile of the sub-sequential limit of $R_{1,2}$. One application of this result is that it simplifies the cavity method as we no longer need a regularizing perturbation to understand the limiting behavior of the overlaps. 
	
	The final motivation of this work is to highlight the hidden role that Panchenko's ultrametricity theorem \cite{PUltra} played in the validity of the cascade conditioning used in \cite{ho_converge,boursier_external_field,chen2026formulafreeenergynonconvex}. In this work, we see that synchronization is a natural property of the Gibbs measure in the limit and is satisfied without one explicitly realizing it or enforcing it through a perturbation. We further explore this connection by providing an interpretation of the cascade conditioning  in Section~\ref{sec:intuition}. Our interpretations are inspired by the limiting picture already known through synchronization \cite{PMult} and ultrametricity \cite{PUltra} and connect the finite $N$ conditioning with its limiting object. We explain how this cascade conditioning approach is a natural construction of ultrametric configurations in the limit by optimally coupling the configurations with the, already ultrametric, uniform Ruelle probability cascades along the depth of the infinite tree. 
	
	\subsection{Setting}
	
	We will consider the multi-species mixed $p$-spin models with Ising spins with an enrichment to match the setting of \cite{boursier_external_field}. We assume no external field for simplicity, but the ideas here should extend to that case. We also expect that these approaches and ideas can generalize to a wide class of spin glass models with Gaussian disorder and sufficient symmetry in the configuration space (see Remark~\ref{rem:bulk}).  
	
	Let $\Sigma_N = \{\pm 1\}^N$ and $\sS$ be a finite set of species labels. For every $N \in \N$, we let $(I_{N,s})_{s \in \sS}$ be a partition of $\{1, \dots, N \}$ and let
	\[
	\lambda_{N,s} := \frac{|I_{N,s}|}{N} \qquad\text{ and } \qquad \vect{\lambda}_N = ( \lambda_{N,s} )_{s \in \sS}
	\]
	denote the corresponding proportions of spins within each species. We assume that $\vect{\lambda}_N$ converges to a vector $\vect{\lambda}$ such that $\lambda_s \in (0,1)$ for all species $s$. Also, define the multi-species overlap 
	\[
	R_{N,s}(\sigma^1,\sigma^2) := \frac{1}{N} \sum_{i \in I_{N,s}} \sigma_i^1 \sigma_i^2, \qquad \vect{R}_N(\sigma^1,\sigma^2) := ( R_{N,s}(\sigma^1, \sigma^2) )_{s \in \sS}.
	\]
	For two replicas $\sigma^\ell, \sigma^{\ell'}$, we often use the shorthand $\vect{R}_{\ell,\ell'} = \vect{R}_{N}(\sigma^\ell,\sigma^{\ell'})$ and $R^s_{\ell,\ell'} = R_{N,s}(\sigma^\ell,\sigma^{\ell'})$ for the species overlap to keep the dependence on $N$ implicit. 
	
	The Hamiltonian has covariance encoded by the function $\xi : \R^{|\sS|} \to \R$ of the form
	\[
	\xi(\vect{x}) = \sum_{k = 1}^\infty \sum_{s_1, \dots, s_k \in \sS} \beta^2_{s_1, \dots, s_k} \prod_{j = 1}^k x_{s_j},
	\] where the coefficients $\beta$ are chosen so that the sum converges absolutely for every $\vect{x}\in \R^{|\sS|}$.
	Let $H_N$ be a centered Gaussian process on $\{ \pm 1 \}^N$ with covariance
	\begin{equation}\label{eq:hamil}
		\E[H_N(\sigma^1) H_N(\sigma^2)] = N \xi( \vect{R}_{1,2} ).   
	\end{equation}
	
	We enrich the multi-species Hamiltonian with a Ruelle probability cascade \cite[Chapter~2.3]{PBook}. We follow the setting of \cite{boursier_external_field}.
	\begin{defn}\label{defn:reg}
		A path is regular if it is continuous and piecewise affine with finitely many pieces and if $q(0) > 0$ and $\dot q > 0$ on every piece. 
	\end{defn}
	Let $q$ be a regular path and $n \geq 1$. For $\alpha \in \N^{n-1}$, let $v_\alpha^{n}$ be the weights of the Ruelle probability cascade whose overlap distribution is uniform on $\{\frac{1}{n},\frac{2}{n},\dots,1\}$. For every regular path $q: [0,1] \to \R$ we can associate it with a sequence
	\[
	0 < q_0 < q_1 < \dots < q_{n} < \infty
	\]
	where $q_k = q(\frac{k}{n})$ for $1 \leq k \leq n$. Let 
	$g(\alpha)$ be the Gaussian process indexed by $\N^{n - 1}$ given by
	\[
	g(\alpha) = q_1^{\frac{1}{2}} g_{\emptyset} + \sum_{k = 1}^{n-1} g_{\alpha | k} ( q_{k+1} - q_{k} )^{1/2}
	\]
	By definition, we have that 
	\[
	\E[ g(\alpha^1) g(\alpha^2) ] =  q( U^{(n)}_{1,2} ) \qquad\text{for}\qquad U^{(n)}_{1,2} = \frac{\alpha^1 \wedge \alpha^2 + 1}{n}.
	\]
	Let $g_i(\alpha)$ be independent copies of $g(\alpha)$ independent of $H_N$. Then, we consider the enriched Hamiltonian given by $$H_N(\sigma,\alpha)= H_N(\sigma) + \sum_{i = 1}^N g_i(\alpha) \sigma_i.$$ This is a centered Gaussian process with covariance \begin{equation}\label{eq:hamil2}   \E[H_N(\sigma^1,\alpha^1)H_N(\sigma^2,\alpha^2)]=N\left(\xi(\vect{R}_{1,2}) + q(U_{1,2})  \sum_{s \in \sS} R_{1,2}^s \right).    
	\end{equation} Let $(\sigma^\ell,\alpha^\ell)_{\ell \geq 1}$ be i.i.d. samples from the corresponding Gibbs measure \begin{equation}\label{eq:Gibbs}
		\dd G_{N}(\sigma, \alpha ) \propto e^{H_N(\sigma,\alpha)} \dd \pP_N(\sigma) v^{(n)}_\alpha ,   
	\end{equation} where $\pP_N(\sigma)$ is the uniform measure on $\Sigma_N$. Furthermore, let $\langle \cdot \rangle$ denote the average with respect to the corresponding product measure $G_{N}^{\otimes \infty}$. In the following, we study the effect of such an enrichment on the corresponding overlap arrays $(\vect{R}_{\ell,\ell'}, U^{(n)}_{\ell,\ell'})_{\ell,\ell'\geq 1}$.
	
	Our main result is the asymptotic synchronization of the $R_N$ and $U^{(n)}$ overlaps along convergent subsequences. More precisely, we get synchronization of the arrays in the limit as $n \to \infty$ followed by $N \to \infty$ along a properly constructed sub-sequence so that the array $\vect{\tilde R} = (\vect{R}_{\ell,\ell'}, U^{(n)}_{\ell,\ell'})_{\ell,\ell'\geq 1}$ converges in distribution. Our theorem takes the $n \to \infty$ limit to get the continuous RPC, and then the $N \to \infty$ to get the synchronized array.
	
	\begin{theo}[$L^2$-Synchronization]\label{thm:main1} Let $q$ be regular. Let $(N_k)$ and $(n_{k,j})$ be such that, for each $k$, the array $(\vect{R}_{\ell,\ell'}, U_{\ell,\ell'}^{(n)})_{\ell \neq \ell'}$ converges in distribution under $\E G_N^{\otimes \infty}$ to $\nu_k$ as $j \to \infty$, and $\nu_k \to \nu$ as $k \to \infty$. We let $(\vect{R}_{\ell,\ell'}, U_{\ell,\ell'})_{\ell \neq \ell'}$ for any array with law $\nu_k$ or $\nu$.
		\begin{enumerate}
			\item  Under $\nu_k$ and $\nu$,  $(U_{\ell,\ell'})$ is a Ruelle probability cascade with continuous uniform overlap distribution.
			\item For $s \in \sS$, let $F_s^{-1}$ be the quantile function of $R_{1,2}^s$ under $\nu$. Then $R_{\ell,\ell'}^s = F_s^{-1}(U_{\ell,\ell'})$ $\nu$ a.s. and for all $\ell \neq \ell'$,$$\lim_{k \to \infty} \lim_{j \to \infty} \E \left\langle(R_{\ell,\ell'}^s-F_s^{-1}(U^{(n)}_{\ell,\ell'}))^2\right\rangle = 0.  $$ 
		\end{enumerate}
	\end{theo}
	
	\begin{rem}
		In this setting, the diagonals of the overlap array converge to a constant by construction, so we also get convergence of the diagonal elements for free. Hence, it is understood in all computations that we are only worrying about the off-diagonal terms. 
	\end{rem}
	\begin{rem}\label{rem:bulk} 
		The proof really only requires the Hamiltonian to be Gaussian with some symmetries such as invariance under permutations or spin flips so this statement can likely hold for a much larger class of models. See the proof of Lemma~\ref{lem:synch} for a detailed explanation.
	\end{rem}
	Synchronization is a natural property of many spin glasses due to the underlying ultrametric structure often enforced through a random regularizing perturbation to imply the Ghirlanda--Guerra identities in the limit on average over the perturbation. This result is surprising because it holds for any covariance function $\xi$, even ones that do not contain the joint moments of the overlaps. In some sense, the structural condition is instead enforced by the enrichment and $q$. This gives intuition on why $q$ must piecewise affine. 
	
	The proof of Theorem~\ref{thm:main1} relies on a very recent cascade conditioning trick to couple configurations with the cascade. We explain the essential idea. Let $(\sigma^\ell,\alpha^\ell)_{\ell \geq 1}$ denote samples from $G_N^{\otimes \infty}$. To synchronize the overlaps, we have to control the conditional variance
	\[
	\E \langle V^s_{\ell,\ell'} \rangle = \E \langle \E\langle (R^s_{\ell,\ell'} - \E \langle R^s_{\ell,\ell'} \given  U_{\ell,\ell'}  \rangle )^2 \given U_{\ell,\ell'}  \rangle \rangle. 
	\]
	The idea is that the set $\{U_{\ell,\ell'} = r \}$ can be understood as  the event that the replicas branch off.  The conditional expected value minimizes the mean squared error, so one might wonder how to attain this estimator. This was done using a branch process in \cite{boursier_external_field}, which can be understood in some sense as an $n \to \infty$ limit of the discrete RPC. This connection was explained in the Appendix of  \cite{boursier_external_field}. The agreement of the cascades at the level of the free energy is in \cite[Section~A.4]{boursier_external_field}. Synchronization gives the agreement at the level of the Gibbs measure. 
	
	\subsection{Background} Multi--species models were first introduced in \cite{mult_og} as a generalization of \cite{SK}. A first step in this direction was the proof of an upper bound on the limit of the free energy  in \cite{mult_og} using Guerra's interpolation method \cite{guerra2003broken}. The crucial matching lower bound was proven in \cite{PMult} using a cavity computation argument. In the lower bound, a new regularizing perturbation was used to show the joint ultrametricity in the limit implies synchronization. However, the upper bound relied crucially on the convexity of $\xi$, which fails for a class of very important models such as the bi--partite spin glass. The approaches for upper and lower bound were extended to spherical multi--species models in \cite{bates_MS} and multi-species Potts models in \cite{ko_cut}.
	
	A brilliant new perspective on free energies motivated by Hamilton--Jacobi equations was proposed in \cite{Mourrat_OG} (see \cite{HJBook} for a textbook introduction) to possibly tackle the non-convex models \cite{Mourrat_noncvx}. These models introduced an enriched spin glass to control the Hamilton--Jacobi error in the limit. The limit free energy for the enriched spin glass was proved in \cite{mourrat_panchenko_sync} using the synchronization developed earlier for multi-species spin glasses. This is also related the synchronization mechanism used in the original proof of the Parisi formula \cite{talagrand2006parisi} alongside Guerra's interpolation. 
	
	Despite the new approach, non-convex $\xi$ remained hard. In spherical models its limit was identified in a subset of non-convex models in \cite{subagTAP,subagTAPII,bates_noncvx,mckenna_noncvx}. One of the first demonstrations of the power of the PDE approach was the conjectured formula for non-convex spin glasses \cite{Mourrat_noncvx} which allowed to identify a critical point condition for the limit free energy \cite{convex_crit,chen_MS} and later a full variational formula \cite{cim_noncvx}, solving one of the big open problems in this field. 
	
	There has been rapid recent progress in these models. Some notable examples being the Hopfield models in \cite{Gaussian_hopfield,barra2026existencethermodynamiclimithopfield}, non-convex spherical models \cite{chen2026freeenergynonconvexmultispecies}, and multi-species spin glasses with external field in \cite{boursier_external_field}. Most remarkably, in \cite{boursier_external_field} synchronization and the Ghirlanda--Guerra identities no longer played a crucial role and the limit free energy was computed using purely the ideas from Hamilton--Jacobi equation. The key ingredient was a conditioning argument on a branch process. Similar ideas involving the cascade conditioning with respect to the local field process \cite{Aukosh_localfield} were introduced earlier in \cite{ho_converge}. The goal of this work is to pick up in another direction where the PDE approach left off and explore the deeper consequences of the cascade conditioning trick.
	
	In this work, we realized that the core idea used to circumvent the Ghirlanda--Guerra identities (GGI) \cite{GGI} and synchronization in the PDE approach is also the core idea that makes it an inherent property of these models. To support the Hamilton--Jacobi approach, this paper adopts a probabilistic perspective tying the picture back to where it started using ultrametricity and revealing the essential role it plays in the validity of these new arguments. Our main contribution is showing that for enriched models, the overlap array synchronizes under the Gibbs measure without a perturbation. 
	
	\subsection{Intuition behind of the proof}\label{sec:intuition} The relatively short length of the paper is deceptive, and is mostly because we use two very complex but important results. The core estimate relies on the control of the conditional variance of one overlap \cite[Theorem~5.1]{boursier_external_field} to give convergence. This involves a complicated cascade conditioning argument along the tree encoded by the enrichment. We do not reproduce that argument here, but we will supplement it with some intuition on the argument and explain why it is intuitive using the limiting ultrametric structure of the Gibbs measures. The other core ingredient is Panchenko's ultrametricity theorem \cite{PUltra} and multi-species synchronization \cite{PMult} to characterize the limit points. 
	
	\subsubsection{Intuition behind the proof}\label{sec:outline_synch}
	
	We will use the results of \cite{boursier_external_field} as a black box in the actual proof, but we will take the time to give our interpretation of why  synchronization is a natural consequence of the conditional overlap concentration results despite never being explicitly identified. 

    This argument first appeared in \cite{aronow2026replicasymmetrybreakingarbitrary} and was used in \cite{ho_converge} to prove the convergence of the free energy. In \cite[Proposition~3.9]{ho_converge}, it was shown that for any species $s,t \in \sS$ with balanced spins the overlaps concentrate
	\begin{equation}\label{eq:concentration_old}
		\E \langle | R_{1,2}^{s} -   R_{1,2}^{t} |^2 \rangle \to 0 .
	\end{equation} Remarkably, no perturbation argument was necessary.
	This statement is in general false unless the overlaps satisfy the following properties:
	\begin{enumerate}
		\item [(a)] The laws of  $R_{1,2}^{s}$ and $R_{1,2}^{t}$ agree, that is $R_{1,2}^{s}\overset{d}{=}R_{1,2}^{t}$ for any species $s,t\in \sS$.
		\item [(b)] The joint distribution of the entire array $(R_{N,s}(\sigma^1,\sigma^2))_{s \in \sS}$ is optimally coupled in the sense that $(R_{N,s}(\sigma^1,\sigma^2))_{s \in \sS} \stackrel{d}{=}  (F_s^{-1}(U))_{s \in \sS}$ for some common source of uniform randomness $U$.
	\end{enumerate}
	Each statement in isolation is believable, the first is obtained through symmetry and the second is attained through synchronization. These mechanisms, although natural, are actually competing against each other. The coupling through synchronization is attained by adding a regularizing perturbation, whose role is to essentially break the symmetries of the model. Yet, for \eqref{eq:concentration_old} to hold, we must have that both symmetry and synchronization hold for all species.
	This seems to suggest that synchronization can be attained without an additional perturbation.
	
	The full power of the cascade conditioning argument  appeared recently in \cite{boursier_external_field,chen2026formulafreeenergynonconvex}. In those works, cascade conditioning was used to show that the error term in the Hamilton--Jacobi equation vanished in the limit without appealing to the Ghirlanda--Guerra identities and synchronization. The previous approaches required introducing a regularizing perturbation to control this term, but this new approach got around this step completely. This seems to suggests that ultrametricity and synchronization are playing a hidden but essential role. 
	
	This series of works are an indication that the enrichment contains enough information on the overlaps to synchronize the overlaps with the enrichment overlap $U$ without perturbation. This was achieved in the past for the generic mixed $p$-spin Hamiltonian, which was rich enough to differentiate in the various temperatures to access the overlap moments. The condition that the cascades order parameters were strictly increasing order parameter is the natural analogue of the dense condition for generic models, because this gives us freedom to perturb the path in any direction and consequently to gain access about all functions of the $U$ overlap. However, we run into a peculiar problem: the derivatives with respect to the order parameter in the direction $\kappa$ give overlaps of the form (see for example \cite[Lemma~3.1]{ho_converge})
	\begin{equation}\label{eq:moments}
		\kappa( U_{1,2} ) \bigg(  \sum_{s \in \sS} R_{1,2}^s \bigg)    
	\end{equation}
	which contain full information about $U$ overlap (which is already known), but almost no information about the powers of the $R$ overlaps which are precisely what we need to control in the GGI. The question is whether the concentration gives a weaker form of synchronization that only holds for the marginals, or whether a stronger mechanism is hidden in the background that couples the entire array.
	
	The answer is the latter. Panchenko's ultrametricity theorem and the  Ghirlanda--Guerra identities allow us to synchronize the entire array just by synchronizing a marginal. This relies crucially on the fact that we are conditioning on the uniform RPC, which is naturally ultrametric. Thus this observation implies that enrichment must optimally couple every sampled pair $(\sigma^\ell,\alpha^\ell)$ at every $N$ to give synchronization of the entire array in the limit. The continuous uniform is the natural object, because all random variables can be generated as the quantile of a continuous uniform random variable. 
	
	Synchronization \cite{PMult,PPotts} is a consequence of ultrametricity \cite{PUltra}. One can interpret the cascade conditioning in \cite{boursier_external_field} as a method to build ultrametricity into the model at finite $N$ instead of merely in the limit through the Ghirlanda--Guerra identities. The idea is simple because we already know the limiting structure through the foundational work of \cite{PUltra}. The enrichment overlap contains the universal structure of all possible limit points of the overlap arrays. By sampling both $(\sigma,\alpha)$ and optimally coupling $\sigma$ with $\alpha$ in the mean-square sense by conditioning, you are forcing $\sigma$ to be as close to a function of $\alpha$ as possible. We know through alternative perturbative methods $\sigma$ is indeed a function of $\alpha$ in the limit, so it is no surprise that the conditional variance control at finite $N$ eventually tends to zero. 
	
	With this intuition in mind, the proof of synchronization is very natural. We can synchronize the overlaps without perturbation using the exact same cascade conditioning in \cite{boursier_external_field}. Rather than recreating the proof from scratch, we can simply apply the bounds as they are relevant in both the PDE and probabilistic approaches. In particular, \cite[Theorem~5.1]{boursier_external_field} contains an estimate on variance of the conditional overlaps,
	\[
	\E\langle \Var( R^s_{\ell,\ell'} \given U_{\ell,\ell'} ) \rangle \to 0 \qquad \text{$s \in \sS$}.
	\]
	
	The connection of the conditional variance to synchronization can now be explained using elementary statistical ideas. The variance of a random variable $X$ is $0$ if it is not random. The conditional variance of $X \given Y$ is $0$, when $X$ is not random conditionally on $Y$, i.e. when $X$ is a function of $Y$.  In particular, we can construct a perfect estimate of $R_{\ell,\ell'}^s$ given $U_{\ell,\ell'}$ in the limit precisely when $R_{\ell,\ell'}^s$ is a function $f_s$ of $U_{\ell,\ell'}$ in the limit. This crucial observation explains why the information in the derivatives \eqref{eq:moments} was able to characterize the entire array.
	
	However, we do not know if $f_s$ is an increasing function, and this is essential for the quantile representation. Convergence of the conditional variance to zero is not enough to give this alone. This is where ultrametricity of the enrichment plays its most crucial role to identify the limit point. 
	
	By exchangeability, we have that all off-diagonal entries of the overlap array $R_{\ell,\ell'}^s$ must be the same function of $U_{\ell,\ell'}$. We know that the array $U$ is ultrametric, so any function of $U$ satisfies the Ghirlanda--Guerra identities. Since the overlaps arrays are Gram de--Finetti arrays, we can use Panchenko's ultrametricity theorem \cite{PUltra} to conclude that the array $R^s$ are ultrametric. This implies that $f_s$ must be a non-decreasing function of $U$, and by uniqueness it is the quantile transform of its marginal distribution. 
	
	This intuitive picture can be made precise by showing that all limit points of the array, satisfy the Ghirlanda--Guerra identities. This allows us to use the full power of synchronization already established in \cite{PMult} to characterize the limit in the usual way. The technical details of these arguments and in particular how to approximate the continuous RPC with a discrete one are given in Section~\ref{sec:synch}.
	
	\subsubsection{Statement of AI use} For pedagogical reasons, AI was not used in the preparation of the initial version of this paper that established the core ideas and proof strategy. This note is the author's interpretation of the cascade conditioning arguments recently used in \cite{ho_converge, chen2026formulafreeenergynonconvex, boursier_external_field}.  Before posting this note, we put an earlier version of the note in Claude Opus 5.5 with the prompt ``You are the referee at a probability journal. Please provide your thoughts on the correctness of these arguments." We did this for two rounds. After the core argument was verified and understandable to the authors, we used AI liberally to edit and polish the final result. 
	
	\subsubsection{Acknowledgments:} The authors sincerely thank Jeanne Boursier, Hong-Bin Chen, Fu-Hsuan Ho, Victor Issa, and Jean-Christophe Mourrat for meeting with the authors and generously explaining their work, sharing their ideas, and providing comments on an earlier version of this work. The open discussion among mathematicians instead of AI inspired this work.
	
	\section{Synchronization} \label{sec:synch}
	
	We first derive $L^2$--synchronization as a consequence of the conditional overlap concentration result of \cite{boursier_external_field} and the multi-species synchronization of \cite{PMult} through the Ghirlanda--Guerra identities \cite{GGI}. In earlier work, the approach to show synchronization was through a perturbation so that the Gibbs measure satisfied the Ghirlanda--Guerra identities either on average or with a carefully chosen sequence of parameters. We will do the same in this work. We start by recalling the Ghirlanda--Guerra identities  for the multi-species model \cite[Theorem~3]{PMult}. 
	\begin{defn}[Ghirlanda--Guerra identities]
		The array $\vect{R} = ( R^1, \dots, R^k)$ satisfies the multi-species Ghirlanda--Guerra identities if for any $m \geq 2$ and any bounded measurable function $f = f( \vect{R}^{\leq m} )$ where $\vect{R}^{\leq m} = ((R^s_{\ell,\ell'})_{s \leq k})_{\ell,\ell' \leq m}$ and any bounded measurable function $\phi:\R^{k}\to\R$ it holds that
		\[
		\E[ f(\vect{R}^{\leq m}) \phi(\vect{R}_{1,m + 1}) ] = \frac{1}{m} \E[ f(\vect{R}^{\leq m})] \E[ \phi(\vect{R}_{1,2}) ] + \frac{1}{m} \sum_{\ell = 2}^m \E[ f(\vect{R}^{\leq m})\phi(\vect{R}_{1,\ell})].
		\]
	\end{defn}
	We start by showing synchronization is equivalent to an array satisfying the GGI. 
	\begin{lem}\label{lem:equivalent}
		Let $U_{\ell,\ell'}$ be a Ruelle probability cascade associated with the continuous uniform distribution. The following are equivalent:
		\begin{enumerate}
			\item[(i)] $( (R^s)_{s\in\sS}, U_{\ell,\ell'} )$ satisfies the multi-species Ghirlanda--Guerra identities;
			\item[(ii)] there exist non-decreasing functions $f_s$ such that
			\[
			\big((R^s_{\ell,\ell'})_{s\in\sS}, U_{\ell,\ell'}\big)_{\ell\neq\ell'}\overset{d}{=}\big((f_s(U_{\ell,\ell'}))_{s\in\sS}, U_{\ell,\ell'} \big)_{\ell\neq\ell'}.
			\]
		\end{enumerate}
	\end{lem}
	\begin{proof}
		(i)$\Rightarrow$(ii): We can essentially treat $U$ as an additional species. Since $((R^s)_{s\in\sS},U)$ satisfies the Ghirlanda--Guerra identities, we have that $$( (R^s_{1,2})_{s\in\sS} , U_{1,2}) \overset{d}{=}(F_{R_{1,2}^{s}}^{-1}(U_{1,2}),U_{1,2})_{s\in \sS}$$ for $U_{1,2} \sim\text{Unif}[0,1]$ (see \cite[Theorem~4]{PMult} or \cite[Section~2]{Mourrat_vec}). Equivalently, it holds that $$R^s_{1,2}=F^{-1}_{R^s_{1,2}}(F_{ R_{1,2}}( R_{1,2})) \quad \text{ a.s. for all }s\in\sS, $$ where $ R_{\ell,\ell'}=\sum_{s\in \sS} R^s_{\ell,\ell'} + U_{\ell,\ell'}$. By weak exchangeability and since there are only countably many pairs $(\ell,\ell')$ and species $s \in \sS$, this implies that  
		\begin{equation}\label{eq: a.s. part} ((R^s_{\ell,\ell'})_{s\in \sS}, U)_{\ell\neq\ell'}=((F^{-1}_{R^s_{1,2}}(F_{ R_{1,2}}( R_{\ell,\ell'})))_{s\in \sS}, F_{ R_{1,2}}( R_{\ell,\ell'}))_{\ell\neq\ell'}  \quad \text{a.s.} 
		\end{equation}
		Notice that $F_{ R_{1,2}}( R_{\ell,\ell'})$ is uniform on $[0,1]$. To conclude from here notice that every element in the array is again weakly exchangeable, positive semi-definite, and satisfies the (usual) Ghirlanda--Guerra identities. In particular, it is ultrametric. Now, notice that the law of $F^{-1}_{ R_{1,2}}(U_{\ell,\ell'})$ is equal to the law of $ R_{1,2}$. Since by \cite[Theorem 2.13]{PBook} the distribution of the overlap matrix under these assumptions is uniquely determined by one overlap value, it must be that \begin{equation}\label{eq: d part}
			( R_{\ell,\ell'}, U_{\ell,\ell'})_{\ell\neq\ell'} \overset{d}{=} (F_{ R_{1,2}}^{-1}(U_{\ell,\ell'}),U_{\ell,\ell'})_{\ell\neq\ell'}.
		\end{equation}
		To finish, combine equations \eqref{eq: a.s. part} and \eqref{eq: d part} and \cite[Lemma 2.2]{Mourrat_vec}.
		
		(ii)$\Rightarrow$(i): Let $f, \phi$ be arbitrary bounded and measurable functions as in the definition of the Ghirlanda--Guerra identities. Define the function $F(u)=((f_s(u) )_{s\in\sS}, u)$. Then, $f\circ F$ and $\phi\circ F$ are again bounded and measurable. Since the uniform Ruelle probability cascade satisfies the Ghirlanda--Guerra identities (see \cite[Section 2.3]{PBook}), we can choose in particular the test functions $f\circ F$ and $\phi\circ F$. Using characterization (ii) then immediately gives that $( (R^s)_{s\in \sS} , U)$ satisfies the Ghirlanda--Guerra identities for arbitrary $f, \phi$.
	\end{proof}
	Our goal is to show that any limiting array satisfies the Ghirlanda--Guerra identities. One approach is the classical approach by showing that the enrichment contains enough information to ensure the validity of the Ghirlanda--Guerra identities by taking derivatives of the enriched free energy and proving an estimate on its rate of concentration \cite[Chapter~3.2]{PBook}. This method has one challenge in the sense that the higher moments of the species overlaps are \emph{not} accessible by differentiating the parameters in our model (see \eqref{eq:moments}). 
	However, this is not an issue because of the crucial realization that the enrichment overlaps $U^{(n)}$ already satisfies the Ghirlanda--Guerra identities. Thus, we essentially get the Ghirlanda--Guerra identities for free if we can show that the species overlaps are a function of the enrichment overlaps $U^{(n)}$. This is accomplished by showing concentration of the conditional variance. 
	
	We define the conditional variance of the array
	\begin{equation}\label{eq:condvar}
		V^s_{\ell,\ell';n} = \E\langle (R^s_{\ell,\ell'} - \E \langle R^s_{\ell,\ell'} \given  U^{(n)}_{\ell,\ell'}  \rangle )^2 \given U^{(n)}_{\ell,\ell'} \rangle.
	\end{equation}
	The conditional expected values are defined in terms of the total expectation $\E \langle \cdot \rangle$, and not only with respect to the Gibbs average $\langle \cdot \rangle$. Our goal will be to control \eqref{eq:condvar}. An approach to this is through  the recently discovered cascade conditioning estimates briefly explained in the introduction. Instead of restating it here, we simply reduce our setting to the branch process setting in \cite[Section~5]{boursier_external_field} and apply those results. Surprisingly, the same estimate that controlled the error in the Hamilton--Jacobi remainder to bypass synchronization is also the same estimate that implies it. 
	\begin{lem}\label{lem:synch}
		If the path is regular, then for every $s \in \sS$ there exists a $C$ so that for all $N \geq 1$, 
		\[
		\lim_{n \to \infty} \E\langle V_{\ell,\ell';n}^s \rangle \leq \frac{C}{N^{1/3}}.
		\]
	\end{lem}   
	
	\begin{proof} Our proof is a reduction of our problem to a form where we can apply \cite[Theorem~5.1]{boursier_external_field} as a black box. The ingredients in the proof of \cite[Theorem~5.1]{boursier_external_field} requires a Gaussian Hamiltonian to get conditional log-concavity in \cite[Lemma~5.6]{boursier_external_field} and permutation invariance in \cite[Proposition~5.4]{boursier_external_field} to remove the conditioning. The exact covariance structure is not needed in the proof of the conditional variance (it is needed to bound the HJ remainder), so our proof will hold for both samples from Gibbs measure and bulk Gibbs measures. 
		
		The main challenge is that the enrichment is of the form of the Ruelle probability cascades, while \cite[Theorem~5.1]{boursier_external_field} uses a branch process. The connections are given in the appendices of \cite{boursier_external_field}. These  play a crucial role in our proof, but were not used in their argument. Our goal is to bound the expected value of the conditional variance $V_{\ell,\ell'}^s$ defined in \eqref{eq:condvar}. 
		
		By \cite[Lemma~A.4]{boursier_external_field} we have that conditioning the configuration on $U^{(n)} = \frac{j}{n-1}$ for $j = 0, \dots, n-1$ is equivalent to conditioning the branch process (see \cite[Section~2]{boursier_external_field}) to split at time $r= \frac{j+1}{n}$ as $n \to \infty$. Furthermore, \cite[Proposition~A.6]{boursier_external_field} shows that the conditional overlaps under the discrete RPC converges, uniformly in depth, as $n \to \infty$ to the conditional overlaps at the corresponding level of the branch. More precisely, \cite[(A.24)]{boursier_external_field} expresses the conditional expected value of the overlap conditioned on the branch at depth $r$ as the limit of the overlap conditioned on the discrete uniform cascade to split at depth $\min(\frac{\lfloor nr \rfloor}{n-1},n-1)$ as $n \to \infty$. We let 
		\[
		\sum_{s \in \sS} \tilde V^s_{\ell,\ell'}(r) :=  \E \langle \| \vect{R}_N - \E \langle \vect{R}_N \rangle_r  \|^2 \rangle_r
		\]
		denote the conditional variance of the $R_{\ell,\ell'}^s$ conditioned to branch at $r$. By the convergence of the conditional overlaps with respect to the uniform cascade to the overlaps conditioned on the branch \cite[Remark~A.7]{boursier_external_field}, we have that 
		\begin{equation}\label{eq:first_estimate}
			\sum_{s \in \sS} \E\langle V^s_{\ell,\ell',n} \rangle = \int  \E \langle \| \vect{R}_N - \E \langle \vect{R}_N \rangle_r  \|^2 \rangle_r \, dr + o_n(1).     
		\end{equation}
		\textit{Step 2:} Applying the conditional variance bound \cite[Theorem~5.1 and Equation 4.9]{boursier_external_field}, implies that
		\[
		\int \E \langle \| \vect{R}_N - \E \langle \vect{R}_N \rangle_r  \|^2 \rangle_r \, dr \leq  \frac{C}{N^{1/3}}
		\]
		Since the conditional variance are non-negative, we have by  \eqref{eq:first_estimate}
		\[
		\E\langle V^s_{\ell,\ell';n} \rangle   \leq \sum_{s \in \sS} \E\langle V^s_{\ell,\ell';n} \rangle  \leq  \int \E \langle \| \vect{R}_N - \E \langle \vect{R}_N \rangle_r  \|^2 \rangle_r \, dr + o_n(1).
		\]
	\end{proof}
	The previous lemma gives us control of one entry of our overlap array in terms of the corresponding entry of $U^{(n)}_{\ell,\ell'}$. To show the Ghirlanda--Guerra identities, we now show that in the limit, the species overlap arrays $R_{\ell,\ell'}$ must be a function of the enrichment array. 
	
	\begin{lem}[Weak Synchronization]\label{eq:weak_sync}
		There exists a sequence of deterministic functions $h_N$ so that for any $\ell,\ell' \geq 1$ and $s \in \sS$,
		\[
		\lim_{N \to \infty} \lim_{n \to \infty} \E \langle (R^s_{\ell,\ell'}  -h^s_{N,n}(U^{(n)}_{\ell,\ell'}) )^2  \rangle \to 0.
		\]
	\end{lem}
	\begin{proof}
		We define $h^s_{N,n}(U^{(n)}_{\ell,\ell'}) = \E\langle R^s_{\ell,\ell'} \given U^{(n)}_{\ell,\ell'} \rangle$. Then by Lemma~\ref{lem:synch},
		\[
		\E \langle (R^s_{\ell,\ell'}  -h_{N,n}^s(U^{(n)}_{\ell,\ell'}) )^2 \rangle = \E \langle (R^s_{\ell,\ell'}  - \E\langle R^s_{\ell,\ell'} \given U^{(n)}_{\ell,\ell'} \rangle )^2 \rangle = \E\langle V_{\ell,\ell'}^s \rangle.
		\]
		The joint array $(\vect{R}_{\ell,\ell'},U^{(n)}_{\ell,\ell'})$ is weakly exchangeable in the $\ell,\ell'$ coordinates by construction, so the same family of function $(h^s_{N,n})_{s \in \sS}$ can be used for off-diagonal entries of the array. 
	\end{proof}
	The coupling can be interpreted as a form of `weak' synchronization in the sense that the species overlaps become deterministic functions of the same enrichment overlap. For many applications, weak synchronization is enough. However, the usual synchronization in spin glasses gives a much richer structure, and in particular that the marginals are ultrametric. We emphasize that the existence a possibly convergent deterministic sequence of functions $h_{N,n}^s$ is in general not enough to conclude synchronization. It is crucial that the functions $h$ is positive and non-deceasing in the limit. 
	
	The key observation is that $U_{\ell,\ell'}^{(n)}$ already satisfies the Ghirlanda--Guerra identities. Thus, the argument in Lemma~\ref{lem:equivalent} part (ii) implies that any function $h_{N,n}^s( U_{\ell,\ell'}^{(n)} )$ also satisfies the Ghirlanda--Guerra identities. Our strategy to show the Ghirlanda--Guerra identities for the overlap arrays $\vect{\tilde{R}} = (\vect{R}_{\ell,\ell'}, U_{\ell,\ell'}^{(n)})$ is now clear. We use Lemma~\ref{eq:weak_sync} to approximate the arrays then use the fact that $(h_{N,n}^s( U_{\ell,\ell'}^{(n)}) , U_{\ell,\ell'}^{(n)})$ satisfies the GGI exactly. 
	
	For any functions $f$ and $\phi$, we define the discrepancy of the array $\vect{\tilde{R}} = ( (R_{\ell,\ell'}^s)_{s \in \sS}, U_{\ell,\ell'} )$ 
	\begin{equation}
		\tilde\Delta(m,f,\phi) = 
		\E\langle f(\vect{\tilde{R}}^{\leq m}) \phi(\vect{\tilde{R}}_{1,m + 1}) \rangle - \frac{1}{m} \E\langle f(\vect{\tilde{R}}^{\leq m})\rangle \E\langle \phi(\vect{\tilde{R}}_{1,2}) \rangle - \frac{1}{m} \sum_{\ell = 2}^m \E\langle f(\vect{\tilde{R}}^{\leq m})\phi(\vect{\tilde{R}}_{1,\ell})\rangle .
	\end{equation}
	
	\begin{lem}\label{lem:GGI}
		For any $m$ and bounded measurable $f,\phi$, we have
		\[
		\lim_{N\to \infty} \lim_{n \to \infty} |	\tilde\Delta(m,f,\phi) | = 0.
		\]
	\end{lem}
	\begin{proof}
		It suffices to show the argument for monomials $f$ and $\phi$, then proceed by approximation. We let $\vect{\tilde{R}}^h = ( (h^s(U^{(n)}_{\ell,\ell'}))_{s \in \sS}, U^{(n)}_{\ell,\ell'} )$ where $h$ is the family of functions $h^s(U^{(n)}_{\ell,\ell'}) = \E\langle R^s_{\ell,\ell'} \given U^{(n)}_{\ell,\ell'} \rangle$. 
		\begin{equation}
		\tilde\Delta^h(m,f,\phi) = 
		\E\langle f(\vect{\tilde{R}}_h^{\leq m}) \phi(\vect{\tilde{R}}^h_{1,m + 1}) \rangle - \frac{1}{m} \E\langle f(\vect{\tilde{R}}_h^{\leq m})\rangle \E\langle \phi(\vect{\tilde{R}}^h_{1,2}) \rangle - \frac{1}{m} \sum_{\ell = 2}^m \E\langle f(\vect{\tilde{R}}_h^{\leq m})\phi(\vect{\tilde{R}}^h_{1,\ell})\rangle .
	\end{equation}
		We first approximate the two discrepancies. Since $f$ and $\phi$ are monomials, there exists an absolute constant $C$ depending on $m,f,\phi$ and all the overlap arrays are bounded, we have
		\[
		| \tilde\Delta^h(m,f,\phi) - \tilde\Delta(m,f,\phi) | \leq \sum_{s \in \sS} C \E \langle |R_{\ell,\ell'}^s - h^s(U^{(n)}_{\ell,\ell'})| \rangle  \leq  C  \sum_{s \in \sS} ( \E\langle V_{\ell,\ell';n}^s \rangle  )^{\frac{1}{2}}
		\]
		since each term is Lipschitz. Since $\vect{\tilde{R}}^h = ( (h^s(U_{\ell,\ell'}^{(n)}))_{s \in \sS}, U_{\ell,\ell'} )$, we have that it satisfies the GGI 
		\[
		\tilde\Delta^h(m,f,\phi)  = 0. 
		\]
		Therefore,
		\[
		\lim_{N\to \infty} \lim_{n \to \infty} | \tilde\Delta(m,f,\phi) |  = \lim_{N\to \infty} \lim_{n \to \infty} | \tilde\Delta^h(m,f,\phi) - \tilde\Delta(m,f,\phi) | = 0. 
		\]
	\end{proof}
	Now our picture is complete. We have that the Ghirlanda--Guerra identities hold in the limit, so we can apply multi-species synchronization \cite{PMult} to characterize all limit points. 
	
	\begin{proof}[Proof of Theorem~\ref{thm:main1}] We first fix the sequence $(N_k, n_{k,j}) $ be a sub-sequence so that for each $k$, $(\vect{R}_{\ell,\ell'}, U^{(n)}_{\ell,\ell'})_{\ell,\ell'\geq 1}$ under the average Gibbs measure to a measure $\nu_k$ as $j \to \infty$ and $\nu_k$ converges to a measure $\nu$.
	\\\\
	\textbf{Part a:} By the Bolthausen--Sznitman invariance \cite[Theorem~2.11]{PBook}, the law of $U^{(n)}_{\ell,\ell'}$ under the average Gibbs measure $\E G_{N,n}$ is a discrete uniform. If the array $\vect{\tilde{R}} = (\vect{R}_{\ell,\ell'}, U_{\ell,\ell'}^{(n)})$ converges as $j \to \infty$ to a measure $\nu_k$, then the distribution of the second coordinate is uniform because it is the limit of discrete uniform random variables. 
	\\\\
	\textbf{Part b:} We have that the Ghirlanda--Guerra identities hold exactly in the limit by Lemma~\ref{lem:GGI}. Therefore, along the convergent sub-sequence as $j \to \infty$, $k \to \infty$, the array $\vect{\tilde{R}} = (\vect{R}_{\ell,\ell'}, U_{\ell,\ell'}^{(n)})$ must also satisfy the Ghirlanda--Guerra identities. Thus by Lemma~\ref{lem:equivalent}, we have that its limit in distribution can be identified by 
		\[
		\big((f_s(U_{\ell,\ell'}))_{s\in\sS}, U_{\ell,\ell'} \big)_{\ell\neq\ell'},
		\]
		where the order parameter of $U$ is the continuous uniform distribution on $[0,1]$. Furthermore, because $U$ is uniform the non-decreasing functions $f_s$ must be the quantile function $F_s^{-1}$ of the limiting distribution of $R^s_{\ell,\ell'}$. 
		
		We can now upgrade to $L^2$ convergence through invariance. We see that
		\begin{equation}\label{eq:L2}
			\E \langle (R^s_{1,2} -  F_s^{-1} (U_{1,2}) )^2 \rangle = 	\E \langle (R^s_{1,2})^2 \rangle  - 2 \E \langle R^s_{1,2}  F_s^{-1} (U_{1,2})  \rangle + \E \langle ( F_s^{-1} (U_{1,2}) )^2 \rangle.
		\end{equation}
		By the Bolthausen--Sznitman invariance \cite[Theorem~2.11]{PBook}, the law of $U^{(n)}_{1,2}$ under $\langle \cdot \rangle$ is uniform. As $j \to \infty$, this converges to a continuous uniform random variable $U_{1,2}$. Furthermore, by synchronization, the joint law of $(R^s_{1,2},U_{1,2})$ is the same as $(F_s^{-1} (U_{1,2}) , U_{1,2})$ since they both generate the same random variable by inverse transform sampling. Thus, in the limit all terms in converge to $\int_0^1 ( F_s^{-1} (u) )^2 \, du$ the second moment of $R^s_{1,2}$, so it vanishes in the limit \eqref{eq:L2}.
	\end{proof}

\end{document}